\documentclass[12pt]{article}
\usepackage{amssymb,latexsym}
\usepackage{amsmath,bbm,amsthm}
\usepackage[active]{srcltx}
\oddsidemargin=-0.1cm \tolerance 9000 \theoremstyle{definition}
\newtheorem{Thm}{Theorem}

\newtheorem{Lem}{Lemma}
\newtheorem{Rem}{Remark}
\newtheorem{Def}{Definition}
\newtheorem{Prop}{Proposition}

\begin{document}

\begin{center}

{{\Large  \textbf{ \sc
Non-local PDEs 
with discrete state-dependent delays: well-posedness in a metric
space}}} 

\vskip7mm

{\large \textsc{Alexander V. Rezounenko }}

\smallskip

Department of Mechanics and Mathematics \\ V.N.Karazin Kharkiv
National  University, 4, Svobody Sqr., Kharkiv, 61077, Ukraine \\
Email: rezounenko@univer.kharkov.ua

\bigskip 

{\large \textsc{ Petr Zagalak }}

\smallskip

Institute of Information Theory and Automation \\ Academy of
Sciences of the Czech Republic, P.O. Box 18, 182\,08 Praha, Czech
Republic \\ Email: zagalak@utia.cas.cz

\end{center}




\bigskip

{\bf Abstract.} Partial differential equations with discrete
(concentrated) state-dependent delays are studied. The existence and
uniqueness of solutions with initial data from a wider linear space
is proven first and then a subset of the space of continuously
differentiable (with respect to an appropriate norm) functions is
used to construct a dynamical system. This subset is an analogue of
\textit{the solution manifold} proposed for ordinary equations in
[H.-O. Walther,  The solution manifold and $C\sp 1$-smoothness for
differential equations
with state-dependent delay, J. Differential Equations, {195}(1),  (2003) 46--65]. 
 The existence of a compact global attractor is proven.

\section{Introduction}

The partial differential equations (PDEs) with delays have attracted
a lot of attention during the last decades as many processes of the
real world (like an automatically controlled furnace, bi-directional
associative memory (BAM) neural networks, reaction-diffusion
processes) can be described by such kind of equations. Studying
these equations is based on the well-developed approaches to the
ordinary differential equations (ODEs) with delays
\cite{Hale_book,Walther_book,Azbelev} and PDEs without delays
\cite{Hadamard-1902,Hadamard-1932,Lions,Lions-Magenes-book}. Under
certain assumptions both types of equations  describe a kind of
dynamical systems that are infinite-dimensional, see
\cite{Babin-Vishik,Temam_book,Chueshov_book} and references therein;
see also \cite{travis_webb,Chueshov-JSM-1992,Cras-1995,NA-1998} and
to the monograph \cite{Wu_book} that are very close to this work.

In many evolution systems arising in applications the presented
delays are frequently \textit{state-dependent} (SDDs). The theory of
such equations, especially the ODEs, is rapidly developping and many
deep results have been obtained up to now (see e.g.
\cite{Walther2002,Walther_JDE-2003,Walther_JDDE-2007,Krisztin-2003,
MalletParet,Walther-JDDE-2010} and also the survey paper
\cite{Hartung-Krisztin-Walther-Wu-2006} for details and references).
The underlying main mathematical difficulty of the theory of PDEs with SDDs
lies in the fact that the functions describing state-dependent delays
are not Lipschitz continuous on the space of continuous functions -
the main space, on which the classical theory of equations with delays is
developed. This implies that the corresponding \textit{initial value problem}
(IVP) is not generally well-posed in the sense of J.~Hadamard \cite{Hadamard-1902,Hadamard-1932}.

The partial differential equations with state-dependent delays were
first studied in
\cite{Rezounenko-Wu-2006} 
(the case of distributed delays, weak solutions),
\cite{Hernandez-2006} (mild solutions, infinite discrete delay), and
\cite{Rezounenko_JMAA-2007} (weak solutions, finite discrete and
distributed delays). An alternative approach to the PDEs with
discrete SDDs is proposed in \cite{Rezounenko_NA-2009}.

This paper is a continuation of the work \cite{Rezounenko_NA-2010} 
and its goal is to study the approach used for ODEs with SDDs
\cite{Walther2002,Walther_JDE-2003,Hartung-Krisztin-Walther-Wu-2006}
in the case of PDEs. The main idea lies in finding a wider space
$Y\supset X$ such that a solution $u : [a,b] \to Y$ be a Lipschitz
function (with respect to a weaker norm of $Y$), and constructing a
dynamical system on a subset of the space  $C([a,b];Y)$. It should
be emphasized that the dynamical system is constructed on a metric
space that is nonlinear. More precisely, the existence and
uniqueness of solutions with initial data from a wider linear space
is proven first and then a subset of the space of continuously
differentiable (with respect to an appropriate norm) functions is
used to construct the aforementioned dynamical system. This subset
is an analogue of \textit{the solution manifold} proposed in
\cite{Walther_JDE-2003}, see also
\cite{Hartung-Krisztin-Walther-Wu-2006}. We use the same class of
non-local in space variables nonlinear PDEs as in
\cite{Rezounenko_NA-2010}.

The paper is organized as follows. The section 2 is devoted to the
formulation of the model.  The proof of the existence and uniqueness
of (strong) solutions for initial functions from a Banach space
forms a main part of the section 3. In the section 4,  an evolution
operator $S_t$ is constructed and its asymptotic properties in
different functional spaces are investigated. The dissipativeness is
obtained in a Banach space, while the existence of a global
attractor is proven on a smaller metric space (the solution
manifold). The choice of this smaller space is different from that
proposed in \cite{Rezounenko_NA-2010}.

\bigskip

\section{The model with discrete state-dependent delay and preliminaries}


Consider the following non-local partial differential equation
with a discrete state-dependent delay $\eta$
\begin{equation}\label{sdd9-1}
\frac{\partial }{\partial t}u(t,x)+Au(t,x)+du(t,x)
=  b\left( [Bu(t-\eta (u_t),\cdot)](x)\right) \equiv \big( F_1(u_t)
\big)(x),\quad x\in \Omega,
\end{equation}

\noindent
where $A$ is a densely-defined self-adjoint positive linear operator
with domain $D(A)\subset L^2(\Omega )$ and compact resolvent, which means
that $A: D(A)\to L^2(\Omega )$ generates an analytic semigroup,
$\Omega \subset \mathbbm{R}^{n_0}$ is  a smooth bounded domain,
$B: L^2(\Omega) \to L^2(\Omega)$ denotes a bounded operator that will be defined later,
$b:\mathbbm{R}\to \mathbbm{R}$ stands for a locally Lipschitz 
map, 
$d \in \mathbbm{R}, d \ge 0$, and the function
$\eta : C([-r,0]; L^2(\Omega)) \to [0,r]\subset \mathbbm{R}_{+}$ denotes
\textit{a state-dependent discrete delay}.
Let $C\equiv C([-r,0];L^2(\Omega))$. Norms defined on $L^2(\Omega)$
and $C$ are denoted by $||\cdot ||$ and $||\cdot ||_C$,
respectively, and $\langle \cdot,\cdot\rangle$ stands for the inner
product in $L^2(\Omega)$. As usually, $u_t\equiv u_t(\theta)\equiv
u(t+\theta)$ for $\theta \in [-r,0].$

\medskip

\begin{Rem} \label{R1}
The operator $B$ may for example be of the following forms (linear operators)
\begin{equation}\label{sdd9-2} [Bv](x)\equiv  \int_\Omega v(y) \widetilde f(x,y) dy ,\quad x\in
\Omega,
\end{equation}
or even simpler
\begin{equation}\label{sdd9-3}
[Bv](x)\equiv \int_\Omega v(y)  f(x-y) \ell (y) dy ,\quad x\in \Omega,
\end{equation}
where $f : \Omega 
\to \mathbbm{R}$ is a smooth function and $\ell \in C^\infty_0 (\Omega)$.
In the last case the nonlinear term in \eqref{sdd9-1} is of the form
\begin{equation}\label{sdd9-4}
\big( F_1(u_t)\big)(x)\equiv  b\left( \int_\Omega u(t-\eta (u_t), y) f(x-y)
\ell (y) dy \right),\quad x\in \Omega .
\end{equation} \hfill $\Box$
\end{Rem}
\bigskip

Consider the equation \eqref{sdd9-1} with the initial condition
\begin{equation}\label{sdd9-ic}
  u|_{[-r,0]}=\varphi 
\end{equation}
and let
\begin{equation}\label{sdd10-33}
H\equiv \left\{ \varphi \in C([-r,0];D(A^{-\frac{1}{2}})) \, | \, \,
\varphi(0)\in D(A^{\frac{1}{2}}) \right\}.
\end{equation}
Let further
$$
||\varphi ||_H\equiv \max_{s\in [-r,0]} ||A^{-\frac{1}{2}}\varphi(s)|| + ||A^{\frac{1}{2}} \varphi(0)||
$$
be a norm defined on the space $H$ and $D(A^\alpha)$ denote the domain of the operator $A^\alpha$. In the sequal the following assumptions will play an important role.

\begin{itemize}
\item[{\bf (H1.$\eta$)}] The discrete delay function $\eta: H\to [0,r]$ is such that
$$
\exists L_\eta > 0, \quad \exists q \ge 0 \quad ~\text{such that}~ \,\, \forall \varphi,\psi\in H \Rightarrow
$$
\begin{equation}\label{sdd10-6}
|\eta(\varphi)-\eta(\psi)|\le L_\eta \left( q ||A^{\frac{1}{2}}
(\varphi(0)-\psi(0))||^2+\int^0_{-r}||A^{-\frac{1}{2}}(\varphi(\theta)
-\psi(\theta))||^2\,d\theta\right)^{\frac{1}{2}}
\end{equation}
\item[{\bf (H.B)}] The following Lipschitz property of the operator $B$ holds.
\begin{equation}\label{sdd9-5}
\exists L_B>0 ~\text{such that}~ \forall u,v\in D(A^{-\frac{1}{2}} ) \Rightarrow
 ||Bu - Bv|| \le L_B  ||A^{-\frac{1}{2}}(u-v)||
\end{equation}
\end{itemize}
\medskip

\begin{Rem} \label{R2}
Under the assumption that for all (almost all) $x\in \Omega
\Rightarrow f(\cdot - x)
\ell (\cdot)\in D(A^{\frac{1}{2}})$ and $u\in 
D(A^{-\frac{1}{2}})$, the term of the form \eqref{sdd9-3} implies that
$$
|\, \langle u,f(\cdot - x)\ell (\cdot) \rangle | \le ||A^{-\frac{1}{2}}u||
\mspace{3mu} ||A^{\frac{1}{2}}f(\cdot - x)\ell (\cdot) ||,
$$
which gives
$$
\left( \int_\Omega \big| \int_\Omega u(y)f(y-x)\ell (y) dy \big|^2\, dx \right)^{\frac{1}{2}}
\le ||A^{-\frac{1}{2}}u||  \left( \int_\Omega || A^{\frac{1}{2}} f(\cdot-x)\ell
(\cdot) ||^2\, dx \right)^{\frac{1}{2}}\mspace{-10mu}.
$$
Hence, the property (H.B) (see \eqref{sdd9-5}) holds with $ L_B \equiv
\left( \int_\Omega || A^{\frac{1}{2}} f(\cdot-x)\ell (\cdot) ||^2\, dx \right)^{\frac{1}{2}}\mspace{-10mu}. $
The same arguments hold (with $ L_B\equiv \left( \int_\Omega ||
A^{\frac{1}{2}} \widetilde f(x,\cdot) ||^2\, dx \right)^{\frac{1}{2}}$) for a more
general term of the form \eqref{sdd9-2}.\hfill $\Box$
\end{Rem}
\medskip


%
%
\medskip

Let now the following space 
\begin{equation}\label{sdd9-8}
{\cal L}\equiv \left\{ \varphi \in C([-r,0];D(A^{-\frac{1}{2}})) \,
|\, \sup\limits_{s\neq t} \left\{ \frac{||A^{-\frac{1}{2}}(\varphi
(s)-\varphi (t)) ||}{|s-t|} \right\} <+\infty;\, \varphi(0)\in
D(A^{\frac{1}{2}})\right\},
\end{equation}
with the natural norm
\begin{equation}\label{sdd9-9}
||\varphi||_{\cal L}\equiv \max_{s\in [-r,0]} ||A^{-\frac{1}{2}}\varphi(s)||
+\sup\limits_{s\neq t} \left\{ \frac{||A^{-\frac{1}{2}}(\varphi (s)-\varphi
 (t)) ||}{|s-t|} \right\} + ||A^{\frac{1}{2}} \varphi(0)||
\end{equation}
be defined.
For any segment $[a,b]\subset \mathbbm{R}$ (c.f. (\ref{sdd9-8}))
and any Lipschitz-on-$[a,b]$ function $\varphi$, let
\begin{equation}\label{sdd10-27}
|||\varphi|||_{[a,b]} \equiv \sup\limits_{} \left\{
\frac{||A^{-\frac{1}{2}}(\varphi (s)-\varphi (t)) ||}{|s-t|} : s\neq t;\,
s,t\in [a,b]\right\}
\end{equation}
denote its Lipschitz constant and let $|||\varphi||| \equiv |||\varphi|||_{[-r,0]}$.
Then the following lemma holds.
\medskip

\begin{Lem} \label{L1}
Let the assumptions (H1.$\eta$)  and (H.B) hold (see
(\ref{sdd10-6}), (\ref{sdd9-5})) and let the function $b:\mathbbm{R}\to \mathbbm{R}$ is
Lipschitz and bounded ($|b(s)| \le M_b$ for all $s\in \mathbbm{R}$).
Then any two functions $\varphi \in {\cal L}, \psi \in H$ (with $H$ and ${\cal L}$ defined
 in (\ref{sdd10-33}) and (\ref{sdd9-8})) the nonlinearity $F$ satisfies
\begin{equation}\label{sdd10-2}
||F_1(\varphi)- F_1(\psi)|| \le L_{F_1} \left[ |||\varphi|||\rule{0pt}{15pt} \right]
\left( q \mspace{3mu} ||A^{1/2}(\varphi(0)-\psi(0))||+ ||A^{-1/2}(\varphi-\psi)||_C\right),
\end{equation}
where
\begin{equation}\label{sdd10-3}
L_{F_1} \left[ \ell \rule{0pt}{1pt} \right]  
\equiv L_b L_B \sqrt{2} \max\left\{ 1; \ell  L_\eta \max \{ 1;
\sqrt{r} \}\right\}
\end{equation}
and $L_{F_1} \left[ \ell \rule{0pt}{1pt} \right]$ is used  in (\ref{sdd10-2}) with
$$ \ell = L_\varphi \equiv |||\varphi||| \equiv
\sup \left\{ \frac{||A^{-1/2}(\varphi (s)-\varphi (t)) ||}{|s-t|} :
s\neq t;\, s,t\in [-r,0]\right\}.
$$
\end{Lem}


\medskip

 \begin{proof}[\textbf{Proof of Lemma 1.}]
Using the Lipschitz property of $b$
 and $B$ (see {\bf (H.B)}), it follows that
\begin{multline*}
|| F_1(\varphi)-F_1(\psi)||^2= \int_\Omega | b([B\varphi](-\eta(\varphi),x))-
 b([B\psi](-\eta(\psi),x)) |^2\, dx \le \\
\le L^2_b \int_\Omega | [B\varphi](-\eta(\varphi),x)-
 [B\psi](-\eta(\psi),x) |^2\, dx = L^2_b || [B\varphi](-\eta(\varphi),\cdot)-
 [B\psi](-\eta(\psi),\cdot) ||^2 \le \\
\le L^2_b L^2_B \mspace{3mu} ||A^{-1/2}\left\{ \varphi
 (-\eta(\varphi))-\psi(-\eta(\psi)) \pm \varphi(-\eta(\psi))\right\} ||^2 \le \\
\le 2L^2_b L^2_B  \left( ||A^{-1/2}\left\{ \varphi
 (-\eta(\varphi))- \varphi(-\eta(\psi))\right\} ||^2 +
 ||A^{-1/2}\left(  \varphi -  \psi\right) ||^2_C \right).
\end{multline*}
Next, $\varphi \in {\cal L}$ implies that there exists
$L_{\varphi}\equiv |||\varphi||| >0,$ (see
(\ref{sdd9-9}),(\ref{sdd10-27}))  such that
\begin{equation}\label{sdd10-18}
||A^{-1/2}(\varphi(s^1)-\varphi(s^2))||\le L_{\varphi}
|s^1-s^2|, \quad \forall s^1,s^2 \in [-r,0].
\end{equation}
Hence, (\ref{sdd10-18}) and (H1.$\eta$) give
$$
|| F_1(\varphi)-F_1(\psi)||^2 \le
$$
\begin{multline*}
\le 2L^2_b L^2_B \bigg[ L^2_{\varphi} L^2_{\eta}
\left( q \mspace{4mu} ||A^{1/2}(\varphi(0)-\psi(0))||^2 + \int^0_{-r} ||A^{-\frac{1}{2}}
(\varphi(\theta)-\psi(\theta))||^2\,d\theta\right)  +  \\
+   || A^{-\frac{1}{2}}\left(  \varphi -  \psi\right) ||^{2}_{C}  \bigg] \le \\
\le 2L^2_b L^2_B  \bigg[ L^2_{\varphi} L^2_{\eta} \left( q \mspace{4mu}
||A^{1/2}(\varphi(0)-\psi(0))||^2 + r \mspace{4mu} ||A^{-\frac{1}{2}}(\varphi-\psi)
||^2_C\,d\theta\right)  + \\
+  ||A^{-\frac{1}{2}}\left( \varphi - \psi\right) ||^2_C  \bigg] \le \\
\le 2L^2_b L^2_B  \max\left\{ 1; L^2_\varphi L^2_\eta \max
\{ 1; r \}\right\}   \left[ q \mspace{4mu} ||A^{1/2}(\varphi(0)-\psi(0))||^2 +
||A^{-\frac{1}{2}} \left( \varphi -  \psi\right) ||^2_C \right ].
\end{multline*}
The last estimate and using the formulas  $\sqrt{\max \{ |a|;|b|\}}
= \max \{ \sqrt{|a|};\sqrt{|b|} \} $ and $\sqrt{a^2+b^2\}} \le |a|+|b|$
give (\ref{sdd10-2}), (\ref{sdd10-3}), which completes the proof.
\end{proof}
\bigskip

\section{The existence and uniqueness of solutions}

\medskip

As in \cite{Rezounenko_NA-2010} we need the following

\medskip

\medskip

\begin{Def}\label{D1}
A vector-function $u(t)\in C([-r,T];D(A^{-1/2}))\cap
C([0,T];D(A^{1/2})) \cap L^2(0,T;D(A))$ with derivative $\dot
u(t)\in L^\infty (0,T;D(A^{-1/2}))$ 
is a (strong) solution to the problem defined by (\ref{sdd9-1}) and
(\ref{sdd9-ic}) on $[0,T]$ if
\begin{itemize}
    \item[(a)] $u(\theta)=\varphi (\theta)$ for $\theta\in [-r,0]$;
    \item[(b)] $\forall v\in L^2(0,T;L^2(\Omega)) $ such that $\dot v\in
L^2(0,T;D(A^{-1})) $ and $v(T)=0 \Rightarrow$
\begin{multline}\label{sdd9-7}
\hskip-10mm -\int^T_0 \langle u(t),\dot v(t)\rangle \, dt  +
\int^T_0 \langle A^{1/2} u(t), A^{1/2} v(t)\rangle \, dt = \\
= \langle\varphi (0),v(0)\rangle + \int^T_0 \langle F_1(u_t)-d\cdot u(t), v(t)\rangle \, dt.
\end{multline}
\end{itemize}
\end{Def}
\medskip


Now we prove the following theorem on the existence and uniqueness of  
solutions.

\medskip

\begin{Thm}\label{T1}
Let the assumptions (H1.$\eta$) and (H.B) hold and let the function
$b:\mathbbm{R}\to \mathbbm{R}$ be Lipschitz and bounded, i.e.
$|b(s)| \le M_b$ for all $s\in \mathbbm{R}$. Let further $\varphi
\in {\cal L}$ be a given initial condition. Then the problem defined
by (\ref{sdd9-1}) and (\ref{sdd9-ic}) has a unique 
solution
on any time interval $[0,T]$ such that $\dot u \in L^2 (0,T;L^2
(\Omega))$.
\end{Thm}
\medskip

\begin{Rem}\label{R3}
Notice that $\varphi$ does not assume
$\varphi\in L^2 ([-r,0];D(A))$. However,  the definition of a strong
solution implies that 
\begin{equation}\label{sdd9-11}
u_t\in L^2 ([-r,0];D(A)), \quad \forall t\ge r.
\end{equation}\hfill $\Box$
\end{Rem}
\medskip

\begin{proof}[\textbf{Proof of Theorem 1.}]
We follow the proof of Theorem~1 in \cite{Rezounenko_NA-2010}.
Notice that the assumption (H1.$\eta$) is slightly more general than
the assumption (H.$\eta$) in \cite{Rezounenko_NA-2010}. This implies
some changes in the proof of the uniqueness of solutions.

Let   $\{ e_k\}^\infty_{k=1}$ denote an orthonormal basis of
$L^2(\Omega)$ such that $Ae_k=\lambda_ke_k$, $0<
\lambda_1<\ldots<\lambda_k\to +\infty$ and consider the Galerkin
approximate solution $u^m=u^m(t,x)=\sum^m_{k=1} g_{k,m}(t) e_k$
of order $m$ such that 
 \begin{equation}\label{sdd9-12}
\left\{ \begin{array}{ll} &\langle \dot u^m+Au^m
+du^m-F_1(u^m_{t}),
e_k\rangle =0,\\
&
\langle u^m(\theta),e_k\rangle=\langle \varphi (\theta) ,
e_k\rangle ,\,\,\forall \theta\in [-r,0] \end{array}\right.
 \end{equation}
$\forall k=1,\ldots,m$, $g_{k,m}\in C^1(0,T;\mathbbm{R})\cap
L^2(-r,T;\mathbbm{R})$ with $\dot g_{k,m}(t)$ absolutely continuous.

The system (\ref{sdd9-12}) is a system of (ordinary) differential
equations in $\mathbbm{R}^m$ with a concentrated (discrete) state-dependent
delay for the unknown vector function $U(t)\equiv (g_{1,m}(t),
\ldots,g_{m,m}(t))$ (for the corresponding theory see
\cite{Walther_JDE-2003,Walther_JDDE-2007} and also a recent review
\cite{Hartung-Krisztin-Walther-Wu-2006}).

The key difference between equations with state-dependent and
state-independent (concentrated) delays is that the first type of
equations is not well-posed in the space of continuous (initial)
functions. To get a well-posed initial value problem, it is better
\cite{Walther_JDE-2003,Walther_JDDE-2007,Hartung-Krisztin-Walther-Wu-2006}
to use a smaller space of Lipschitz continuous functions or even
a smaller subspace of $C^1([-r,0];\mathbbm{R}^m)$.

The condition $\varphi \in {\cal L}$ implies that the function
$U(\cdot)|_{[-r,0]}\equiv P_m \varphi (\cdot) $, which defines
initial data, is Lipschitz continuous as a function from $[-r,0]$ to
$\mathbbm{R}^m$. Here $P_m$ is the orthogonal projection onto the
subspace $span\, \{ e_1,\ldots, e_m\}\subset L^2(\Omega)$.  Hence,
we can apply the theory of ODEs with discrete state-dependent delay
(see e.g.
\cite{Hartung-Krisztin-Walther-Wu-2006}) to get the local
existence and uniqueness of solutions to (\ref{sdd9-12}).

Next, we will get an a priory estimate to prove the continuation of
solutions $u^m$ to (\ref{sdd9-12}) on any time interval $[0,T]$
and then use it for the proof (by the method of compactness,
see \cite{Lions}) of the existence of strong solutions to
(\ref{sdd9-1}) and (\ref{sdd9-ic}).  
To that end, multiply the first equation in (\ref{sdd9-12}) by $\lambda_k
g_{k,m}$ and sum for $k=1,\ldots,m$ to get
$$\frac{1}{2} \frac{d}{dt} ||A^{1/2}u^m(t)||^2 + ||Au^m(t)||^2 +
d\cdot ||A^{1/2}u^m(t)||^2= \langle P_m F(u^m_t),Au^m(t)\rangle \le
$$
$$\le \frac{1}{2} ||P_m F(u^m_t)||^2 + \frac{1}{2} ||Au^m(t)||^2.
$$
As the function $b$ is bounded, $||F(u^m_t)||^2\le M^2_b |\Omega|$
(here $|\Omega|\equiv \int_\Omega 1 \, dx$), which gives

\begin{equation}\label{sdd9-13}
\frac{d}{dt} ||A^{1/2}u^m(t)||^2 + ||Au^m(t)||^2 \le M^2_b
|\Omega|.
\end{equation}
Integrating (\ref{sdd9-13}) with respect to $t$ and using the relationships
$\varphi (0)\in D(A^{1/2}),$ 
$u^m(0)=P_m\varphi (0)\in D(A^{1/2})$,
$||A^{1/2}u^m(0)||=||A^{1/2}P_m\varphi (0)|| \le ||A^{1/2}\varphi
(0)||$, we get an a priory estimate
\begin{equation}\label{sdd9-14}
||A^{1/2}u^m(t)||^2 + \int^t_0 ||Au^m(\tau)||^2 \, d\tau \le
||A^{1/2}\varphi (0)||^2 + M^2_b |\Omega|\mspace{3mu} T, \quad \forall m,
\forall t\in [0,T].
\end{equation}
\medskip

The above relationship (\ref{sdd9-14}) means that
$$
\{ u^m \}^\infty_{m=1} \hbox{ is a bounded set in }  L^\infty(0,T;D(A^{1/2}))\cap L^2(0,T;D(A)).
$$
Using this fact and (\ref{sdd9-12}), it follows that
$$
\{ \dot u^m \}^\infty_{m=1} \hbox{ is a bounded set in }  L^\infty(0,T;D(A^{-1/2}))\cap L^2(0,T;L^2(\Omega)).
$$
Hence, the family $\{ (u^m;\dot u^m ) \}^\infty_{m=1}$ is a bounded set in
\begin{multline}\label{sdd9-15}
Z_1\equiv \left( L^\infty(0,T;D(A^{1/2}))\cap L^2(0,T;D(A))\right)\times \\
\times\left( L^\infty(0,T;D(A^{-1/2}))\cap L^2(0,T;L^2(\Omega)) \right).
\end{multline}
Therefore, there exist a subsequence $\{ (u^k;\dot u^k ) \}$ and an
element  $(u;\dot u )\in Z_1$ such that
\begin{equation}\label{sdd9-16}
\{ (u^k;\dot u^k ) \} \hbox{ *-weak converges to } (u;\dot u )
\hbox{ in } Z_1.
\end{equation}
The proof that any *-weak limit is a strong solution is standard.
 To prove 
the property $u(t) \in C([0,T];D(A^{{1/2}}))$, we use the well-known
(see also \cite[thm. 1.3.1]{Lions-Magenes-book}) 
\medskip

\begin{Prop}\label{P2} (Proposition 1.2 in \cite{showalter}).
Let  $V$ denote a dense Banach space that is continuously embedded in
a Hilbert space $X$ and let  $X=X^*$ so that $V\hookrightarrow
X\hookrightarrow V^*$. Then the Banach space $W_p(0,T)\equiv \{ u\in
L^p(0,T;V) : \dot u\in L^q(0,T;V^*)\}$ (here $p^{-1}+q^{-1}=1$) is
contained in $C([0,T];X).$
\end{Prop}
\medskip

\noindent In our case $X=D(A^{{1/2}})$, $V=D(A)$, $V^*=L^2(\Omega)$, $p=q=1/2$.
\medskip

Now we prove the uniqueness of 
 solutions. Using the fact that
$\varphi \in {\cal L}$, the definition~1 of a 
solution $v$,
and $\dot v(t) \in L^\infty (0,T; D(A^{-1/2}))$ (see
(\ref{sdd9-16})), it follows that for any such a solution $v$ and any
$T>0$, there exists $L_{v,T}>0,$ such that
\begin{equation}\label{sdd9-18}
||A^{-1/2}(v(s^1)-v(s^2))||\le L_{v,T}  |s^1-s^2|, \quad
\forall s^1,s^2 \in [-r,T].
\end{equation}
In the light of (\ref{sdd10-27}), let $L_{v,T}\equiv |||v|||_{[-r,T]}$.
\smallskip

Consider any two 
solutions $u$ and $v$ of (\ref{sdd9-1}),
(\ref{sdd9-ic}) (not necessarily with the same initial function).
The standard variation-of-constants formula $ u(t)=e^{-At} u(0) +
\int^t_0 e^{-A(t-\tau)} F(u_\tau)\, d\tau$
and the estimate  $|| A^\alpha e^{-tA}|| \le \left(\frac{\alpha}
{t}\right)^\alpha e^{-\alpha}$ (see e.g. \cite[(1.17),
p.84]{Chueshov_book}) give
$$
 ||A^{1/2}(u(t) - v(t))|| \le
 e^{-\lambda_1t}||A^{1/2}(u(0) - v(0))|| + \int^t_0 || A^{1/2}
e^{-A(t-\tau)}|| \mspace{5mu} ||F(u_\tau) -F(v_\tau) || \, d\tau \le
$$
\begin{equation}\label{sdd10-20}
\le e^{-\lambda_1t}||A^{1/2}(u(0) - v(0))|| + \int^t_0
\left(\frac{1/2}{t-\tau}\right)^{1/2} e^{-1/2} \mspace{3mu}
||F(u_\tau) -F(v_\tau) || \, d\tau,
 \end{equation}
as $|| A^{1/2} e^{-A(t-\tau)}|| \le \left(\frac{1/2}{t-\tau}\right)^{1/2} e^{-1/2}$
and similarly,
$$
||A^{-1/2}(u_t - v_t)||_C \le ||A^{-1/2}(u_0 - v_0)||_C + \int^t_0
||F(u_\tau) -F(v_\tau) ||  \, d\tau.
$$
\medskip

The last estimate and (\ref{sdd10-20})  give
(just the case when $q=1$ is shown for the purpose of clarity)
\begin{multline}\label{sdd9-20b}
||A^{1/2}(u(t) - v(t))|| + ||A^{-1/2}(u_t - v_t)||_C \le
e^{-\lambda_1t}||A^{1/2}(u(0) - v(0))|| + \\
+  ||A^{-1/2}(u_0 - v_0)||_C + \int^t_0 \left\{ 1+ (2e (t-\tau))^{-1/2}\right\}
 ||F(u_\tau) -F(v_\tau) ||  \, d\tau.
\end{multline}

It follows, from Lemma~\ref{L1}, that
\begin{equation}\label{sdd10-17}
||F(u_t)- F(v_t) || \le L_{F_1,v,T}  \left( 
 q\mspace{3mu}      
 ||A^{1/2}(u(t) - v(t))||
+ ||A^{-1/2}(u_t-v_t)||_C \right),
\end{equation}
where $L_{F_1,v,T}$ is defined in the same way as $L_{F_1}$ in (\ref{sdd10-3}),
just with $\ell = L_{v,T}$ instead of $L_{\varphi}$ - see (\ref{sdd10-3})
and (\ref{sdd9-18}).
\begin{equation}\label{sdd10-19}
L_{F_1,v,T} \equiv L_{F_1} 
 \left[ L_{v,T}\rule{0pt}{15pt} \right]  
\equiv  L_b L_B \sqrt{2} \max\left\{ 1; L_{v,T} L_\eta \max \{ 1;
\sqrt{r} \}\right\}.
\end{equation}



It should be emphasized  how the Lipschitz constant
$L_{v,T}\equiv |||v|||_{[-r,T]}$ of a strong solution $v$ is taken
into account in (\ref{sdd10-19}) (see (\ref{sdd9-18}) and
(\ref{sdd10-27})).
\medskip

Let
$$
g(t)\equiv 
 ||A^{1/2}(u(t) - v(t))|| + ||A^{-1/2}(u_t - v_t)||_C.
$$
Then the relationships (\ref{sdd9-20b}) and (\ref{sdd10-17}) lead to
the following estimate
$$ g(t) \le g(0) + \int^t_0 \left\{ 1+ 
 (2e (t-\tau))^{-1/2}\right\} L_{F_1,v,T}\cdot g(\tau)\, d\tau
$$

\begin{Lem}[\textbf{Gronwall}] \label{L2}
Let $u,\alpha \in C[a,b], \beta (t)\ge 0, \beta$ is integrable on $[a,b]$ and
$$
u(t) \le \alpha(t) + \int^t_a \beta (\tau) u(\tau)\, d\tau, \quad a\le t\le b
$$
Then
$$u(t) \le \alpha(t) + \int^t_a \beta (\tau) \alpha(\tau)
\exp\left\{ \int^t_\tau \beta (s) \, ds\right\}\, d\tau, \quad
a\le t\le b
$$
Moreover, if $\alpha$ is non-decreasing, then
$$
u(t) \le \alpha(t)  \exp\left\{ \int^t_a  \beta (s) \, ds\right\},~~   a\le t\le b.
$$
\end{Lem}
\bigskip

It follows, from the above lemma and equality
$\int^t_0 (t-\tau)^{-1/2}d\tau = 2t^{1/2}$, that
\begin{multline*}
g(t) \le g(0) \mspace{3mu} \exp\left\{ L_{F_1,v,T}  \int^t_0
\left\{ 1+ (2e (t-s))^{-1/2}\right\} \, ds\right\} \le \\
\le g(0) \mspace{3mu} \exp\left\{ L_{F_1,v,T}  \left( t+ \sqrt{\frac{2t}{e}}\right) \,\right\},
\end{multline*}
which implies, $\forall t\in [0,T]$, that
\begin{multline}\label{sdd10-20a}
||A^{1/2}(u(t) - v(t))|| + ||A^{-1/2}(u_t - v_t)||_C \le \\
\le E_{F_1,v,T} \left(  ||A^{1/2}(u(0) - v(0))|| +  ||A^{-1/2}(u_0 - v_0)||_C\right),
\end{multline}
where
\begin{equation}\label{sdd10-21}
E_{F_1,v,T} \equiv \exp\left\{ L_{F_1,v,T} \cdot  \left( T+
\sqrt{\frac{2T}{e}}\right) \,\right\},
\end{equation}
see (\ref{sdd10-19}) for the definition of  $L_{F_1,v,T} \equiv L_{F_1}
 \left[ L_{v,T}\rule{0pt}{10pt} \right]$.
This proves the uniqueness of the 
solution to (\ref{sdd9-1}) and (\ref{sdd9-ic}), and completes the
proof of theorem~1.
\end{proof}
\bigskip

\section{Asymptotic properties of solutions}
\medskip

This section is devoted to studies of the asymptotic behavior of solutions in
different functional spaces.
 We define first (in a standard way) the evolution
semigroup $S_t : {\cal L} \to {\cal L}$ (the space ${\cal L}$ is
defined in (\ref{sdd9-8})) by the formula
\begin{equation}\label{sdd10-26}
S_t \varphi \equiv u_t,\quad t\ge 0,
\end{equation}
where $u(t)$ is a unique 
solution to the problem (\ref{sdd9-1}) and (\ref{sdd9-ic}) (see
definition~1).
\medskip

The estimate (\ref{sdd10-20a}) means the continuity of the evolution
operator $S_t$  in the norm of the space $H$ (see (\ref{sdd10-33})), i.e.
\begin{equation}\label{sdd10-25}
||S_t \varphi - S_t \psi ||_H \le E_{F_1,v,T} \mspace{3mu} || \varphi - \psi
||_H \hbox{ for all } t\in [0,T].
\end{equation}

The aim now is to get a more precise estimate, e.g. the continuity
of $S_t$  in the norm of the space ${\cal L}$ (see (\ref{sdd9-8}),
(\ref{sdd9-9})). Consider the definition of the Galerkin approximate
solution (see (\ref{sdd9-12})). It gives
\begin{multline*}
||A^{-1/2}(\dot u^m(t) - \dot v^m(t)) || \le ||A^{1/2}(u^m(t) -  v^m(t)) ||
+ d ||A^{-1/2}(u^m(t) - \\
- v^m(t)) || + ||F_1(u^m_t)-F_1(v^m_t) ||
\end{multline*}
and Lemma~\ref{L1}  implies
\begin{multline*}
||A^{-1/2}(\dot u^m(t) - \dot v^m(t)) || \le (1+d+L_{F1}) \{ ||A^{1/2}(u^m(t) -  v^m(t))|| + \\
+ ||A^{-1/2}(u^m_t -  v^m_t) ||_C \}.
\end{multline*}

An analogous estimate for a 
solution to the problem
(\ref{sdd9-1}) and (\ref{sdd9-ic}), can be obtained from
(\ref{sdd9-16}) and  the following

\smallskip

\begin{Prop}\cite[Theorem~9]{yosida} \label{P3}
Let $X$ be a Banach space. Then any *-weak convergent sequence
 $\{ w_k\}^\infty_{n=1}\in X^{*}$  *-weak converges to an element
 $w_\infty\in X^{*}$ and $\Vert w_\infty\Vert_X \le\liminf_{n\to\infty}
\Vert w_n\Vert_X.$
\end{Prop}
\medskip

\noindent More precisely,
\begin{multline*}
\hbox{ess sup}_{t\in [0,T]} ||A^{-1/2}(\dot u(t) - \dot v(t)) ||
\le (1+d+L_{F1}) \sup_{t\in [0,T]} \{ ||A^{1/2}(u(t) - v(t))||  + \\
+ ||A^{-1/2}(u_t -  v_t) ||_C \}
\end{multline*}
The last estimate and relationship (\ref{sdd10-20a}) imply
\begin{multline} \label{sdd10-20b}
\hbox{ess sup}_{t\in [0,T]} ||A^{-1/2}(\dot u(t) - \dot v(t)) || \le \\
\le (1+d+L_{F1})  E_{F_1,v,T}  \left(  ||A^{1/2}(u(0) -
v(0))|| +  ||A^{-1/2}(u_0 - v_0)||_C\right)
\end{multline}
Hence, see (\ref{sdd10-27}),
$$
||| u - v|||_{[0,T]} \le (1+d+L_{F1}) E_{F_1,v,T} \left(  ||A^{1/2}(u(0) - v(0))||
+  ||A^{-1/2}(u_0- v_0)||_C\right)
$$
From that and (\ref{sdd10-20a}), it follows that
\begin{equation}\label{sdd10-28}
||u_t -v_t||_{{\cal L}} \le (2+d+L_{F1})  E_{F_1,v,T}\mspace{3mu}
||u_0 -v_0||_{{\cal L}},~~ \forall t\in [0,T],
\end{equation}
which finally means that for any $T\ge 0$ there exists a
constant $C_T>0$ such that $\forall t\in [0,T]$ it gives
\begin{equation}\label{sdd10-29}
||u_t -v_t||_{{\cal L}} = ||S_t \varphi -S_t \psi ||_{{\cal L}} \le
C_T  ||\varphi -\psi||_{{\cal L}},~~\forall \varphi, \psi \in {\cal L}
\end{equation}
The last inequality means the continuity of the evolution operator $S_t$ in the
norm of the space ${\cal L}$ (see (\ref{sdd10-26}) and compare with
(\ref{sdd10-25})).
\medskip

\begin{Rem}\label{R4}
It should be noted that the evolution operator and, more generally, the time-shift is not
a (strongly) continuous mapping in the norm of the space ${\cal L}$ (see (\ref{sdd9-8})).
This can be illustrated by the following simple (scalar) example.

\noindent Consider the space
$$ {\cal L}ip\, ([-r,T];\mathbbm{R}) \equiv \left\{ v : [-r,T]\to \mathbbm{R} :
\sup \left\{ \frac{|v(s)-v(t)|}{|s-t|}, s\neq t; s,t\in [-r,T]\right\}
 <\infty \right\}
$$
and analogously define the space $ {\cal L}ip\, ([-r,0];\mathbbm{R})$ with the
natural norm
$$ ||v||_{{\cal L}ip} \equiv \max_{\theta\in [-r,0]} |v(\theta)| +
\sup \left\{ \frac{|v(s)-v(t)|}{|s-t|}, s\neq t;\, s,t\in
[-r,0]\right\}.
$$
The (strong) continuity of the time-shift  means that
\begin{equation}\label{sdd10-34}
\forall v\in {\cal L}ip\, ([-r,T];\mathbbm{R}) \text{ and } \forall
t \in [0,T] \Longrightarrow \lim\limits_{h\to 0} ||v_{t+h} -
v_t||_{{\cal L}ip}=0.
\end{equation}
Obviously, when $t=0$ one considers $h\to 0^{+}$, while for $t=T$, the case $h\to 0^{-}$
should be investigated.

To prove the claim, we must show that  (\ref{sdd10-34}) does not hold, i.e.
\begin{equation}\label{sdd10-35}
\exists v\in {\cal L}ip\, ([-r,T;\mathbbm{R}) \text{ and } \exists
t_0 \in [0,T] \quad \text{for which} \quad \lim\limits_{h\to 0}
||v_{t_0+h} - v_{t_0}||_{{\cal L}ip} \neq 0.
\end{equation}
Thus, consider the case $t_0=0, h\to 0^{+}$ and the function
$$ v(t)\equiv %
\begin{cases}
  0, & t\in [-r,0] \\
  t, & t\in (0,T] \\
\end{cases}.
$$
It can be seen that
$v_{t_0}=v_0\equiv 0$ and
$$
v_{t_0+h} =v_{t_0+h}(\theta) =
\begin{cases}
  0, & \theta\in [-r,-h] \\
  h+\theta, & \theta\in (-h,0]
\end{cases}.
$$

\noindent Hence, $||v_{t_0+h} - v_{t_0}||_{{\cal L}ip} = ||v_{t_0+h}
||_{{\cal L}ip} = h+1 $ and finally $\lim\limits_{h\to 0+}
||v_{t_0+h} - v_{t_0}||_{{\cal L}ip} = \lim\limits_{h\to 0+} (h+1)=1
\neq 0$, which means that (\ref{sdd10-34}) does not hold.
In the space ${\cal L}$, we would proceed analogously. \hfill $\Box$
\end{Rem}
\medskip

\begin{Rem}\label{R5}
In the same way as in the previous remark one
can show that  the time-shift is {\bf not} a (strongly) continuous
mapping in the topology  of $L^\infty (-r,0).$ One could consider
the function $\widetilde v(t)\equiv
\begin{cases}
  0, & t\in [-r,0] \\
  1, & t\in (0,T] \\
\end{cases}$
and  $t_0=0$ to show that $\lim\limits_{h\to 0+} ||\widetilde v_{h}
- \widetilde v_{0}||_{L^\infty (-r,0)} = 1 \neq 0$. By the way,
$\widetilde v = \frac{d}{dt} v,$ where, as usually, the
time-derivative is understood in the sense of distributions. \hfill
$\Box$
\end{Rem}
\medskip

The above remarks show that despite of the existence and uniqueness
of solutions in the space ${\cal L} $ and even strong continuity of
the evolution operator $S_t$ in the norm of ${\cal L}$ (see
(\ref{sdd10-29})), the pair $(S_t;{\cal L}) $ does  {\it not} form a
{\it dynamical system} since $S_t$ is not strongly continuous as a
mapping of time variable.

The methods, developed for ordinary delay equations in
\cite{Walther_JDE-2003} 
suggest to restrict our considerations to a smaller subset of the
space of Lipschitz functions. In this paper we follow this
suggestion and consider the evolution operator  $S_t$  on the
following subset of~${\cal L} $
\smallskip

$ X \equiv \left\{ \quad \varphi \in C^1([-r,0];D(A^{-1/2}))  \qquad
\hbox{ such that } \right.$
\begin{equation}\label{sdd10-36} \left. \varphi
(0)\in D(A^{1/2})) \hbox{ and } \dot \varphi (0) + A \varphi (0) +
d\varphi (0)= F_1(\varphi) \quad \right\} \subset {\cal L}.
\end{equation}
Here the equality $\dot \varphi (0) + A \varphi (0) + d\varphi (0)=
F_1(\varphi)$ is understood as an equality in $D(A^{-1/2})$.
\medskip

\begin{Rem}\label{R6}
The set $X$ is an analogue of the solution
manifold introduced in \cite{Walther_JDE-2003} for the case of ODEs
with state-dependent delays. \hfill $\Box$
\end{Rem}
\medskip

To show that the set $X$ is invariant under the evolution operator
$S_t$, we first have to establish an additional smoothness property
of the solutions of problem (\ref{sdd9-1}), (\ref{sdd9-ic}).
\medskip

\begin{Lem}\label{L3}
For any $\varphi \in C^1([-r,0];D(A^{-1/2}))$ such that $\varphi
(0)\in D(A^{1/2})),$ the 
solution to (\ref{sdd9-1}),
(\ref{sdd9-ic}) (which is given by Theorem~\ref{T1}) has the property (c.f.
Proposition~\ref{P2} and Theorem~\ref{T1})
\begin{equation}\label{sdd10-37}
\dot u \in C([0,T];D(A^{-1/2})), \quad \forall T>0.
\end{equation}
\end{Lem}

\medskip

\begin{Rem}\label{R7}
We do not assume $\varphi (0)\in D(A),$ just
$\varphi (0)\in D(A^{1/2}),$ so we cannot directly use \cite[Theorem
3.5, p.114]{Pazy-1983}.
\end{Rem}
\medskip

\begin{proof}[\textbf{Proof of Lemma~\ref{L3}.}]
By Proposition~1 and Theorem~1, for any
$\varphi \in C^1([-r,0];D(A^{-1/2}))$ such that $\varphi (0)\in
D(A^{1/2}))$, there exists a unique solution $u(t)\in
C([-r,T];D(A^{-1/2}))\cap C([0,T];D(A^{1/2}))$. This property and
Lemma~\ref{L1} then imply the continuity of the function
\begin{equation}\label{sdd10-38}
p(t)\equiv F_1(u_t) \in C([0,T]; L^2(\Omega)). 
\end{equation}
Consider the following auxiliary {\it linear system without
delay}
\begin{equation}\label{sdd10-39}
\left\{ 
\begin{array}{l}
  \dot v(t) +A v(t) + d v(t)=p(t),\quad t\ge 0, \\
  v(0)=\varphi (0)\in D(A^{1/2}) \\
\end{array}
\right.
\end{equation}
In the same way as in (\ref{sdd9-12}), the Galerkin approximate solution
$v^m=v^m(t,x)=\sum^m_{k=1} g_{k,m}(t) e_k$ of order $m$ to (\ref{sdd10-39})
can be defined such that 
 \begin{equation}\label{sdd10-40}
\left\{ \begin{array}{ll}
\langle \dot v^m+Av^m +dv^m-p(t),
e_k\rangle =0,\quad t\ge 0,\\
\langle v^m(0),e_k\rangle=\langle \varphi (0) , e_k\rangle, \quad
\forall k=1,\ldots,m.
 \end{array}\right.
 \end{equation}
where $g_{k,m}\in C^1(0,T;\mathbbm{R})\cap L^2(-r,T;\mathbbm{R})$
and $\dot g_{k,m}(t)$ is absolutely continuous.

The difference between approximate solutions $u^m$ and
$v^m$ lies in that $v^m$ are solutions just to {\it linear} system
(\ref{sdd10-40}). So, for any two approximate solutions $v^n$ and
$v^{m}$ (solutions to (\ref{sdd10-40}) of different orders $n$ and
$m$), one has $g_{k,n} (t)\equiv g_{k,m} (t)$, which is denoted by
$g_{k} (t)$.

Multiply (\ref{sdd10-40}) by $\lambda_k g_k$ and sum for
$k=n+1,...,n+p$ ($p$ is any positive integer) to get
\begin{multline*}
\langle \dot v^{n+p}(t)- \dot v^{n}(t),  A (v^{n+p}(t)-  v^{n}(t)) \rangle +
|| A (v^{n+p}(t)-  v^{n}(t)) ||^2 + \\
+ d \langle v^{n+p}(t)-  v^{n}(t), A (v^{n+p}(t)-  v^{n}(t))
\rangle = \langle (P_{n+p} - P_n)p(t), A (v^{n+p}(t)-v^{n}(t))\rangle
\end{multline*}
It should be recalled that, see the proof of Theorem~\ref{T1}, $P_m$
is the orthogonal projection onto the subspace $span\, \{ e_1,\ldots,
e_m\}\subset L^2(\Omega)$. Hence,
\begin{multline*}
\frac{1}{2} \frac{d}{dt} || A^{\frac{1}{2}}(v^{n+p}(t)-
v^{n}(t)) ||^2 + || A (v^{n+p}(t)- v^{n}(t)) ||^2 + d || A^{\frac{1}
{2}}(v^{n+p}(t)- v^{n}(t)) ||^2  \le \\
\le || (P_{n+p} - P_n)p(t)|| \cdot || A (v^{n+p}(t)-
v^{n}(t))|| \le \frac{1}{2} || (P_{n+p} - P_n)p(t)||^2 + \\
+ \frac{1}{2} || A (v^{n+p}(t)- v^{n}(t))||^2
\end{multline*}
which gives
$$
\frac{d}{dt} || A^{\frac{1}{2}}(v^{n+p}(t)- v^{n}(t)) ||^2 + ||
A(v^{n+p}(t)- v^{n}(t)) ||^2 \le || (P_{n+p} - P_n)p(t)||^2.
$$
Integrating the last estimate results ($\forall t\in [0,T]$) in
\begin{multline*}
|| A^{\frac{1}{2}}(v^{n+p}(t)-v^{n}(t)) ||^2 +
\int^t_0 || A(v^{n+p}(\tau)- v^{n}(\tau)) ||^2\, d\tau \le \\
\le || A^{\frac{1}{2}}(v^{n+p}(0)- v^{n}(0)) ||^2  +
\int^t_0 || (P_{n+p} - P_n)p(\tau)||^2 \, d\tau \le \\
\le || (P_{n+p}- P_{n}) A^{\frac{1}{2}} \varphi(0)) ||^2  + \int^T_0
|| (P_{n+p} - P_n)p(\tau)||^2 \, d\tau.
\end{multline*}

Summing up, the above estimate, the fact that $\varphi (0)\in D(A^{1/2}))$,  the strong
convergence $||I -P_n|| \to 0$ for $n\to \infty$, and (\ref{sdd10-38})
imply that
\begin{equation}\label{sdd10-41}
\textit{the sequence } \{ v^n \}^\infty_{n=1}\textit{  is a Cauchy
sequence in } C([0,T];D(A^{\frac{1}{2}})).
\end{equation}

Now our goal is to show that the sequence $\{ \dot v^n\}^\infty_{n=1} $
is a Cauchy sequence in $C([0,T];D(A^{-1/2}))$. So, multiply first (\ref{sdd10-40})
by $\lambda^{-\frac{1}{2}}_k$
to get $\lambda^{-\frac{1}{2}}_k \dot g_k (t) = - \lambda^{\frac{1}{2}}_k  g_k (t)
- d \lambda^{-\frac{1}{2}}_k g_k (t) + \langle \lambda^{-\frac{1}{2}}_k
 p(t),e_k\rangle$. This gives
 $\lambda^{-1}_k (\dot g_k (t))^2 \le 3 \lambda_k ( g_k (t))^2 + 3d^2 \lambda^{-1}_k (g_k (t))^2
 + 3| \langle \lambda^{-\frac{1}{2}}_k  p(t),e_k\rangle |^2$.
The sum for $k=n+1,...,n+p$ reads
 \begin{multline*}
|| A^{-\frac{1}{2}}(\dot v^{n+p}(t)-\dot v^{n}(t)) ||^2 \le 3 || A^{\frac{1}{2}}
(v^{n+p}(t)- v^{n}(t)) ||^2+ \\
+ 3d^2 || A^{-\frac{1}{2}}(v^{n+p}(t)- v^{n}(t)) ||^2 + \frac{3}
{\lambda_{n+1}} || (P_{n+p} - P_n)p(t)||^2 \le \\
\le 3\left( 1+ \frac{d^2}{\lambda^2_{n+1}}\right) || A^{\frac{1}{2}}
(v^{n+p}(t)- v^{n}(t)) ||^2 +\frac{3}{\lambda_{n+1}} ||I-P_n||^2
||p(t)||^2.
\end{multline*}
The last estimation together with (\ref{sdd10-41}) give that
\begin{equation}\label{sdd10-42}
\hbox{ the sequence } \{ \dot v^n \}^\infty_{n=1} \hbox{ is a Cauchy
sequence in } C([0,T];D(A^{-\frac{1}{2}})).
\end{equation}
Thus, there exists a unique solution $v(t)$ ($v\equiv
\lim_{n\to\infty} v^n$) to the linear system (\ref{sdd10-39}),
which satisfies $v\in C([0,T];D(A^{\frac{1}{2}}))$ and $\dot v\in
C([0,T];D(A^{-\frac{1}{2}}))$.

On the other hand, the nonlinear delay system (\ref{sdd9-1}),
(\ref{sdd9-ic}) with the initial function $\varphi$ has also a
unique solution. From the construction of $p(t)$ (see
(\ref{sdd10-38})), it follows that $u(t)\equiv v(t)$ for all $t\in
[0,T]$, which gives (\ref{sdd10-37}) and completes the proof of
Lemma~\ref{L3}.
\end{proof}
\medskip

Lemma~\ref{L3} particularly shows that the set $X$, defined by
(\ref{sdd10-36}), is invariant under the evolution operator $S_t$
(see (\ref{sdd10-26})). This fact allows to define an evolution
operator (denoted again by $S_t$) $S_t : X \to X$ in the same way as
in (\ref{sdd10-26}). Now, if  the natural norm
$$
||\varphi||_X\equiv \max_{s\in [-r,0]} ||A^{-1/2}\varphi(s)||
+\max_{s\in [-r,0]} ||A^{-1/2}\dot \varphi(s)|| + ||A^{1/2}
\varphi(0)||
$$
on $X$ is taken into account, then Theorem~\ref{T1}, Lemma~\ref{L3},
and Proposition~\ref{P2} give the continuity of $S_t$ with respect
to $t$  in the norm of $X$. Hence, $(S_t;X)$ defines a dynamical
system.
\bigskip

Now we will pay attention to the long-time asymptotic behavior of
the constructed evolution semigroup $S_t : X \to X$.
\begin{Thm}\label{T2}
Using the above notation and under the assumptions of Theorem~\ref{T1},
the dynamical system $(S_t,X)$ is dissipative. If,
in addition, $q=0$ in {\bf (H1.$\eta$)}, then 
$(S_t,X)$ possesses a compact global attractor ${\cal A}$, which is
a bounded set in the space $C^1([-r,0];D(A^{-1/2}))\cap
C([-r,0];D(A^{\alpha}))$, $\alpha \in ({\frac{1}{2}},1)$.
\end{Thm}
\medskip

\begin{proof}[\textbf{Proof of Theorem~\ref{T2}.}]
It will be shown first that $(S_t,X)$ is a
dissipative dynamical system. To that end, the below proposition
is needed. 
\begin{Prop}\label{P4} \cite[Lemma~1]{Rezounenko_NA-2010}
Let all the assumptions of Theorem~\ref{T1} hold and let
$\alpha \in ({\frac{1}{2}},1)$. Then there exists a bounded subset
${\cal B}V_\alpha$
of the space $C^1([-r,0];D(A^{-\frac{1}{2}}))\cap C([-r,0];D(A^{\alpha}))$,
which absorbs any strong solution to the problem
(\ref{sdd9-1}) and (\ref{sdd9-ic}) for any initial function
$\varphi \in {\cal L}$.
\end{Prop}

Second, to apply the classical theorem on the existence of a global
attractor (see, for example
\cite{Babin-Vishik,Temam_book,Chueshov_book}), we show that
$(S_t,X)$ is asymptotically compact. Consider therefore any solution
$u(t)$ to the problem (\ref{sdd9-1}) and (\ref{sdd9-ic}) with
$\varphi \in {\cal B}V_\alpha$ as an initial
function. We will show that for any
$\delta >r>0$ and any $T> \delta$ the set ${\cal U} \equiv
\{ u_t=S_t\varphi \, | \, \varphi \in {\cal B}V_\alpha, \, t\in
[\delta, T] \}$ is relatively compact in $X$.

Recall that the set ${\cal B}V_\alpha$ is a ball in
$C^1([-r,0];D(A^{-1/2}))\cap C([-r,0];D(A^{\alpha}))$ (for more
details see \cite{Rezounenko_NA-2010}) and notice that, by
Corollary~4 from \cite{Simon}, the set  ${\cal B}V_\alpha $ is
relatively compact in $C([-r,0];D(A^{-1/2}))$ (see also
\cite[lemma~1]{Simon}).
It remains to show that $\{ \dot u(t) \, | \, \varphi \in {\cal
B}V_\alpha, \, t\in [\delta-r, T] \}$ is 
equi-continuous  in $C([\delta-r,T];D(A^{-1/2}))$.
\medskip


\begin{Prop}\label{P5} \cite[Corollary 4.3.3 and Theorem 4.3.5]{Pazy-1983}.
Let $A$ be an infinitesimal generator of
an analytic semigroup $\{ T(t)\}_{t\ge 0}$. If $f\in L^1((0,T);Y) $
is locally H$\ddot{o}$lder continuous on $(0,T]$, then for every
$x\in Y$ the initial value problem
$$ \dot u (t)= A u(t) + f(t), t > 0;\quad u(0)=x$$
has a unique solution $u$. If $f\in C^\theta ([0,T];Y)$, then for
every $\delta>0,$ $Au \in C^\theta ([\delta,T];Y)$ and $\dot u \in
C^\theta ([\delta,T];Y)$.
\end{Prop}
\bigskip

\noindent Here $C^\theta ([0,T];Y)$ denotes the family of all H$\ddot{o}$lder
continuous functions on $[0,T]$ with the exponent $\theta \in (0,1)$.
In this case, $Y=L^2(\Omega)$.
\bigskip

In order to apply Proposition~\ref{P5} to our case, we have to show that
$p(t)=F_1(u_t) \in C^\theta ([\delta-r,T];L^2(\Omega))$ (c.f.
(\ref{sdd10-38}),(\ref{sdd10-39})). Therefore, consider $t\in [\delta-r,T]$ and
\begin{multline*}
||p(t+h)-p(t)||= ||F_1(u_{t+h})-F_1(u_t) || \le \\
\le L_{F_1}[\ell_{{\cal B}V_\alpha}\mspace{-2mu}] \max_{s\in [-r,0]}
||A^{-1/2}(u(t+h+s) - u(t+s)) || \le L_{F_1}\left[\ell_{{\cal
B}V_\alpha}\mspace{-2mu}\right] \mspace{3mu}\ell_{{\cal B}V_\alpha}\mspace{3mu} |h|
\end{multline*}
where $L_{F_1}[\ell_{{\cal B}V_\alpha}\mspace{-2mu}]$ is the constant defined in
Lemma~\ref{L1} with $\ell_{{\cal B}V_\alpha}$ such that
$|||\psi|||\le \ell_{{\cal B}V_\alpha}$ $\forall \psi\in {\cal B}V_\alpha$
(the existence of such $\ell_{{\cal B}V_\alpha}$ follows from
Proposition~\ref{P3}). Here, $q=0$ is used.

The last inequality shows that $p : [\delta-r,T] \to L^2(\Omega)$
is Lipschitz continuous, which is the situation to which Proposition~\ref{P4}
can be applied.
It should also be noted that the family $\{ p(t)\}$, for all
initial $\varphi \in {\cal B}V_\alpha$, is uniformly Lipschitz,
i.e. all the Lipschitz constants are lower or equal to $L\equiv
L_{F_1}[\ell_{{\cal B}V_\alpha}]\cdot \ell_{{\cal B}V_\alpha}$. Then
by Proposition~\ref{P3}, it is guaranteed (see the proof)
that the family $\{ \dot u(t) \, | \, \varphi \in
{\cal B}V_\alpha, \, t\in [\delta-r, T] \}$ is uniformly H$\ddot o$lder
continuous, and thus 
equi-continuous  in $C([\delta-r,T];L^2(\Omega))$.
\begin{Prop}\label{P6} \cite[lemma~1]{Simon}
Let $B$ be a Banach space. A set $F$ of $C([0,T];B)$ is relatively
compact if and only if
\begin{itemize}
\item[(i)] $F(t)\equiv \{ f(t) : f\in F\}$ is relatively compact in $B, ~0<t<T$,
\item[(ii)] $F$ is uniformly equicontinuous, i.e. $\forall\varepsilon>0,
\exists \eta$ such that $||f(t_2)-f(t_1)||_B \le \varepsilon,
\forall f\in F, \forall 0\le t_1\le t_2 \le T$ such that $|t_2-t_1 |
\le \eta$
\end{itemize}
\end{Prop}
\bigskip

Applying Proposition~\ref{P6} completes the proof of
Theorem~\ref{T2}.
\end{proof}
\bigskip

As an application we can consider the diffusive Nicholson
blowflies equation (see e.g. \cite{So-Yang}) with state-dependent
delays, i.e. the equation (\ref{sdd9-1}) where $-A$ is the Laplace
operator with the Dirichlet boundary conditions, $\Omega\subset
\mathbbm{R}^{n_0}$ is a bounded domain with a smooth boundary,
the nonlinear (birth) function $b$ is given by $b(w)=p\cdot we^{-w}$.
The function $b$ is bounded
, so for any delay function $\eta$ satisfying $(H1.\eta)$, the
conditions of Theorem~\ref{T1} and Theorem~\ref{T2} are satisfied.
As a result, we conclude
that the initial value problem (\ref{sdd9-1}) and (\ref{sdd9-ic}) is
well-posed in $X$ and the dynamical system $(S_t,X)$ has a global
attractor (Theorem~\ref{T2}).

\end{document}